\documentclass[11pt,reqno]{amsart}
\usepackage{amsmath, amssymb, amsthm}
\usepackage{url}
\usepackage[breaklinks]{hyperref}
\renewcommand{\(}{\left\(}
\renewcommand{\)}{\right\)}
\renewcommand{\[}{\left\[}
\renewcommand{\]}{\right\]}

\numberwithin{equation}{section}
 \theoremstyle{plain}
\newtheorem{theorem}{Theorem}[section]
\newtheorem{lemma}[theorem]{Lemma}
\newtheorem{remark}[]{Remark}

\newtheorem{definition}{Definition}[section]

   \makeatletter
\def\proof{\@ifnextchar[{\@oproof}{\@nproof}}
\def\@oproof[#1][#2]{\trivlist\item[\hskip\labelsep\textit{#2 Proof of\
#1.}~]\ignorespaces}
\def\@nproof{\trivlist\item[\hskip\labelsep\textit{Proof.}~]\ignorespaces}

\makeatother

\newenvironment{proof-alt}
   {\vskip 0.15in \par\noindent{\it Proof of Proposition \ref{HSS}.}\hskip 0.5em\ignorespaces}
   {\hfill $\Box$\par\medskip}

\begin{document}
\title[Congruences for 2-colored partitions]{Congruences for the difference of even and odd number of parts of the cubic and some analogous partition functions}

\author{Nayandeep Deka Baruah}
\address{Nayandeep Deka Baruah, Department of Mathematical Sciences, Tezpur University, Napaam, Assam, 784028, India.}
\email{nayan@tezu.ernet.in, nayandeeptezu@gmail.com}

\author{Abhishek Sarma}
\address{Abhishek Sarma, Department of Mathematical Sciences, Tezpur University, Napaam, Assam, 784028, India.}
\email{abhitezu002@gmail.com}


\maketitle

\begin{abstract}
Partitions wherein the even parts appear in two different colours are known as cubic partitions. Recently, Merca introduced and studied the function $A(n)$, which is defined as the difference between the number of cubic partitions of $n$ into an even number of parts and the number of cubic partitions of $n$ into an odd number of parts. In particular, using Smoot's \textsf{RaduRK} Mathematica package, Merca proved the following congruences by finding the exact generating functions of the respective functions. For all $n\ge0$, \begin{align*}A(9n+5)\equiv 0\pmod 3,\\
A(27n+26)\equiv 0\pmod 3. \end{align*} By using generating function manipulations and dissections, da Silva and Sellers proved these congruences and two infinite families of congruences modulo 3 arising from these congruences. In this paper, by employing Ramanujan's theta function identities, we present simplified formulas of the generating functions from which proofs of the congruences of Merca as well as those of  da Silva and Sellers follow quite naturally. We also study analogous partition functions wherein multiples of $k$ appear in two different colours, where $k\in\{3,5,7,23\}$.\end{abstract}

\section{Introduction}
A partition $\lambda=(\lambda_1, \lambda_2, \ldots, \lambda_k)$ of a positive integer $n$ is a non-increasing sequence of positive integers, $\lambda_1\geq \lambda_2\geq \cdots \geq \lambda_k$ such that $\lambda_1+\lambda_2+\cdots+\lambda_k=n$. Here each $\lambda_{i}$ is called a part of the partition. The partition function $p(n)$ counts the number of partitions of $n$. With the convention $p(0)=1$, the generating function of $p(n)$ due to Euler is given by
\begin{align*}
\sum_{n=0}^{\infty} p(n) q^n = \dfrac{1}{(q;q)_{\infty}},
\end{align*}
where for complex numbers $a$ and $q$ with $|q|<1$, the customary $q$-products are defined by
\begin{align*}(a;q)_0&:=1, ~ (a;q)_n:=\prod_{k=0}^{n-1}(1-aq^k)\quad \textup{for a positive integer}~ n,\\
(a;q)_\infty&:=\lim_{n\to\infty}(a;q)_n.
\end{align*}

The cubic partition function $a(n)$ counts the number of partitions of a positive integer $n$ in which even parts can appear in two colors. This function was introduced by  Hei-Chi Chan [\cite{zbMATH05718299}, \cite{zbMATH05758244}] in 2010 and connected to the so-called Ramanujan's cubic continued fraction. With the convention of $a(0)=1$, the generating function of $a(n)$ is given by
\begin{align*}
	\sum_{n=0}^{\infty} a(n)q^n = \frac{1}{(q;q)_{\infty}(q^2;q^2)_{\infty}}.
\end{align*}

Motivated by Ramanujan's so-called ``most beautiful identity" 
\begin{align*}
	\sum_{n\geq0} p(5n+4)q^n = 5\frac{(q^5,q^5)^5_{\infty}}{(q;q)^6_{\infty}},
\end{align*}
Chan \cite{zbMATH05718299} proved the following analogous identity:
\begin{align*}
	\sum_{n\geq0} a(3n+2)q^n = 3\frac{(q^3,q^3)^3_{\infty}(q^6,q^6)^3_{\infty}}{(q,q)^4_{\infty}(q^2,q^2)^4_{\infty}}.
\end{align*}
This clearly gives 
\begin{align*}
	a(3n+2) = 0 \pmod3,
\end{align*}
which is analogous to
\begin{align*}
	p(5n+4) = 0 \pmod5,
\end{align*}
one of the three famous congruences for the partition function discovered by Ramanujan.

Recently, Merca \cite[Definition 1]{zbMATH07576463} defined the following functions.
\begin{definition}
	For a positive integer $n$, let
	\begin{enumerate} 
		\item $a_{e}(n)$ be the number of partitions of $n$ into an even number of parts in which the even parts can appear in two colours.
		\item  $a_{o}(n)$ be the number of partitions of $n$ into an odd number of parts in which the even parts can appear in two colours.
		\item $A(n)$ = $a_{e}(n)$ - $a_{o}(n)$.
	\end{enumerate}
\end{definition}
By considering \begin{align*}
	F(z,q) = \prod_{n=1}^{\infty} \frac{1}{(1-zq^n)(1-zq^{2n})},
\end{align*}
Merca found that
\begin{align}\label{genAn}
	\sum_{n=0}^{\infty} A(n) q^n = F(-1,q)=\frac{1}{(-q;q)_{\infty}(-q^2;q^2)_{\infty}} = (q;q^2)_{\infty}(q^2;q^4)_{\infty},
\end{align}
where the last equality arises from Euler's identity  \begin{align*}
	(-q,q)_\infty=\frac{1}{(q,q^2)_\infty}.
\end{align*}

Let $t_{1}+2t_{2}+\cdots+nt_{n}=n$ be a partition. From this partition, we can derive partitions of $n$ into exactly $t_{1}+t_{2}+\cdots+t_{n}$ parts where even parts can appear in two colours. The number of such partitions is
\begin{align*}
	(1+t_{2})(1+t_{4})\cdots(1+t_{2[n/2]}).
\end{align*} 
Using this idea, Merca \cite[Theorem 1.1]{zbMATH07576463} deduced the following theorem:
\begin{theorem}\label{th1}
	Let $n$ be a positive integer. Then
	\begin{align*}
		a_{e}(n)\pm a_{o}(n) = \sum_{t_{1}+2t_{2}+\cdots+nt_{n}=n}^{} (\pm1)^{t_{1}+t_{2}+\cdots+t_{n}} (1+t_{2})(1+t_{4})\cdots(1+t_{2[n/2]}).
	\end{align*}
\end{theorem}
Clearly, $a_{e}(n) + a_{o}(n)$ represents the total number of cubic partitions of $n$.

With the aid of Smoot's \textsf{RaduRK} Mathematica package \cite{smoot} which is based on an algorithm developed by Radu \cite{radu1}, Merca \cite[Theorem 6.1]{zbMATH07576463} proved that
\begin{align}\label{merca-gen}
\sum_{n=0}^\infty &A(9n+5)q^n\notag\\&=-3q\dfrac{(q^2;q^2)_\infty^3(q^3;q^3)_\infty(q^{12};q^{12})_\infty^6}{(q;q)_\infty^2(q^4;q^4)_\infty^3(q^{6};q^{6})_\infty^5}+9q^2\dfrac{(q^2;q^2)_\infty^5(q^3;q^3)_\infty(q^{12};q^{12})_\infty^{10}}{(q;q)_\infty^2(q^4;q^4)_\infty^7(q^{6};q^{6})_\infty^7}\notag\\
&\quad+3q^3\dfrac{(q^2;q^2)_\infty^7(q^3;q^3)_\infty(q^{12};q^{12})_\infty^{14}}{(q;q)_\infty^2(q^4;q^4)_\infty^{11}(q^{6};q^{6})_\infty^9}-9q^4\dfrac{(q^2;q^2)_\infty^9(q^3;q^3)_\infty(q^{12};q^{12})_\infty^{18}}{(q;q)_\infty^2(q^4;q^4)_\infty^{15}(q^{6};q^{6})_\infty^{11}}
\end{align}
and a similar expression for the generating function of $A(27n+26)$ having twelve terms \cite[Theorem 6.1]{zbMATH07576463}. From these expressions of the generating functions, Merca readily found the following two Ramanujan-like congruences.
\begin{theorem}\cite[Theorem 1.10]{zbMATH07576463}\label{thm1.2}
	For all $n \geq0$,
	\begin{align}\label{cong-A9n5}
	 A(9n+5) &\equiv 0 \pmod3,\\
	 \label{cong-A27n26}A(27n+26)& \equiv 0 \pmod3.
	\end{align}
\end{theorem}
Using classical generating function manipulations and dissections, da Silva and Sellers \cite{da2022elementary} reproved the above congruences. They also proved an additional congruence and couple of infinite families of congruences modulo 3 as stated in the following theorem.
\begin{theorem}\cite[Theorems 3.1--3.3]{da2022elementary} \label{thm1.3}For all $j \geq 0$ and $n \geq 0$,
\begin{align}\label{Silva-Sellers-Cong1}
			A(3n+1) &\equiv  A(27n+8)\pmod3,\\
		\label{Silva-Sellers-Cong2}A\left(9^{j+1}n + \frac{39\cdot9^j + 1}{8}\right) &\equiv 0 \pmod3,\\
		\label{Silva-Sellers-Cong3}A\left(3\cdot9^{j+1}n + \dfrac{23\cdot9^{j+1} + 1}{8}\right) &\equiv 0 \pmod3.
	\end{align}
\end{theorem}

The first purpose of this paper is to employ Ramanujan's theta function identities in finding simplified formulas of the generating functions from which proofs of Theorem \ref{thm1.2} and Theorem \ref{thm1.3} follow quite naturally. To state our results, we now define Ramanujan's theta functions. Ramanujan's general theta function $f(a,b)$ is defined by 
\begin{align*}
	f(a,b) = \sum_{k=-\infty}^{\infty} a^{\frac{k(k+1)}{2}}b^{\frac{k(k-1)}{2}}, \quad |ab| < 1.
\end{align*}
Three special cases of $f(a,b)$ are:
\begin{align}
	\varphi(q) &:= f(q,q) = (-q;q^2)^2_\infty(q^2;q^2)_\infty,\label{varphi}\\
	\psi(q) &:= f(q,q^3) = \frac{(q^2;q^2)_\infty}{(q;q^2)_\infty},\label{psi}\\
	f(-q) &:= f(-q,-q^2) = (q;q)_\infty,\notag
\end{align}
where the product representations follow from Jacobi's triple product identity \cite[p. 35, Entry 19]{bcb3},
 \begin{align}\label{jtpi}
	f(a,b) = (-a;,ab)_\infty,(-b;,ab)_\infty(ab;,ab)_\infty.
\end{align}
We also define 
\begin{align}\label{chi}
\chi(q):=(-q;q^2)_\infty.
\end{align}
For brevity, we also set $f_k:=(q^k;q^k)_\infty$.

Now we state our results on $A(n)$.
\begin{theorem}\label{thm1.4}We have
	\begin{align}\label{gen-A3n}
		\sum_{n=0}^{\infty} A(3n) q^{n} &= \frac{\chi(-q)\psi(q^3)\varphi(-q^3)}{\varphi(-q)\varphi(q^2)},\\
	\label{gen-A3n1}\sum_{n=0}^{\infty} A(3n+1) q^{n} &= -\frac{\psi^2(q^3)\chi^2(-q)}{\psi(-q)\psi(q^2)},\\
		\label{gen-A3n2}\sum_{n=0}^{\infty} A(3n+2) q^{n} &= -\frac{\psi(-q^3)\psi(q^6)}{\psi^2(q^2)},\\
		\label{gen-A9n5}\sum_{n=1}^{\infty} A(9n+5)q^n &= -3q\dfrac{(q;q)_\infty(q^{2};q^{2})_\infty^{4}(q^{12};q^{12})_\infty^{6}}{(q^{4};q^{4})_\infty^{11}}.\\\intertext{Furthermore, for all $n\geq0$}
		\label{cong-3n2-27n17}A(3n+2)&\equiv-A(27n+17)\pmod3\\
				\label{cong-81n44}A(81n+44) &\equiv 0 \pmod3.
\end{align}
\end{theorem}
Note that \eqref{gen-A9n5} is a much simplified form of \eqref{merca-gen} and the congruence \eqref{cong-A9n5} also readily follows from \eqref{gen-A9n5}.

The second purpose of this paper is to study the partition function $A_k(n)$ defined below, where $A_2(n)=A(n)$. 
\begin{definition}
	For a positive integer $n$, let
	\begin{enumerate}
		\item $a_e^k(n)$ be the number of partitions of $n$ into an even number of parts in which the parts that are multiples of $k$ can appear in two colours.
		\item  $a_o^k(n)$ be the number of partitions of $n$ into an odd number of parts in which the parts that are multiples of $k$ can appear in two colours.
		\item $A_k(n)$ = $a_e^k(n)$ - $a_o^k(n)$.
	\end{enumerate}
\end{definition} 

For example, if $k$= 3 and $n$ = $4$, then  
\begin{align*} a_e^3(4) &= 4, \textup{the relevant partitions being} ~3_r+1, 3_b+1, 2+2, 1+1+1+1;\\
a_o^3(4) &= 2, \textup{the relevant partition being}~ 4~ \textup{and}~ 2 + 1 + 1;\end{align*}
and hence, $A_3(4)=4-2 = 2$, where the subscripts $r$ and $b$ depict the two colors of the respective part.

Note that, Theorem \ref{th1} can be extended to 2-colored partitions where multiples of $k$ can appear in 2-colors. Using the same idea used by Merca \cite[Theorem 1.1]{zbMATH07576463}, one can easily deduce the following theorem.
\begin{theorem}
	For any positive integer $n$, we have
	\begin{align*}
		a_e^k(n)\pm a_o^k(n) = \sum_{t_{1}+2t_{2}+\cdots+nt_{n}=n}^{} (\pm1)^{t_{1}+t_{2}+\cdots+t_{n}} (1+t_{k})(1+t_{2k})\cdots(1+t_{k[n/k]}).
	\end{align*}
\end{theorem}
Clearly, $a_{e}^k(n) + a_{o}^k(n)$ represents the total number of partitions of $n$ where $k$ can appear in two colors.

By considering the function 
\begin{align*}
	G(z,q) = \prod_{n=1}^{\infty} \frac{1}{(1-zq^n)(1-zq^{kn})},
\end{align*}
we can easily see that the generating function of $A_k(n)$ is $G(-1,q)$; that is,  
\begin{align*}
	\sum_{n=0}^{\infty} A_{k}(n) q^n = \dfrac{1}{(-q;q)_{\infty}(-q^{k};q^{k})_{\infty}},
\end{align*}
which by Euler's identity $(-q;q)_\infty=1/(q;q^2)_\infty$ and \eqref{chi} may be recast as
\begin{align}\label{gen-Akn}
	\sum_{n=0}^{\infty} A_{k}(n) q^n = (q;q^2)_\infty(q^k;q^{2k})_\infty=\chi(-q)\chi(-q^k).
\end{align}

In the following theorems we state our results on $A_k(n)$ for $k\in\{3,5,7,23\}$.

\begin{theorem}\label{thm1.6}
For all $n \geq 0$, 
	\begin{align}\label{cong-A3-4nr}
		A_{3}(4n+r) \equiv 0 \pmod2,\quad where\quad r \in \{2, 3\}.
	\end{align}	
\end{theorem}

\begin{theorem}\label{thm1.7}
For all $n \geq 0$, 
	\begin{align}\label{cong-A5-10nr}
		A_{5}(10n+r) &\equiv 0 \pmod{2},\quad where\quad r \in \{2, 6\},\\		\label{cong-A5-25nr}A_{5}(25n+r) &\equiv 0 \pmod{5},\quad where\quad r \in \{14, 19, 24\}.
	\end{align}	
\end{theorem}

\begin{theorem}\label{thm1.8}
For all $n\geq0$ and $j\geq0$,\begin{align}
		\label{cong-A7-2n1-8n3}
		A_7(2n+1)&\equiv A_7(8n+3) \pmod2,\\
		\label{cong-A7-inf8n7}
		A_7\left(2^{2j+3}n+\dfrac{10\cdot2^{2j+1}+1}{3}\right)&\equiv0\pmod2,\\
		\label{cong-A7-16nr}A_{7}(16n+r) &\equiv 0 \pmod{2},\quad where\quad r \in \{9,13\},
		\end{align}
		for all $\alpha\geq0$,
		\begin{align}
\label{cong-A7-14nr-inf}
A_7\left(2\cdot7^{\alpha}(7n+r)+\dfrac{2\cdot7^{\alpha}+1}{3}\right) &\equiv 0 \pmod2,\quad where\quad r \in \{3,4,6\},
\end{align}
			and if $p>3$ is a prime such that $\left(\frac{-7}{p}\right)= -1$, then for all $\alpha\geq0$,
		\begin{align}\label{cong-A7-pnr}A_7\left(2p^{2\alpha+1}(pn+r)+ \dfrac{2p^{2\alpha+2}+1}{3}\right) \equiv 0 \pmod2,		
	\end{align}
	where $r \in \{1,2,\ldots,p-1\}$.
	\end{theorem}

	\begin{theorem}\label{thm1.9}
	For all $n\geq0$,
	\begin{align}\label{cong-A23-23nr}
		A_{23}\left(2\cdot23^\alpha (23n+r)+2\cdot23^\alpha+1\right)\equiv 0\pmod2,\end{align}
		where $r\in\{4,6,9,10,13,14,16,18,19,20,21\}$,
		
		\noindent and if $p>3$ is a prime such that $\left(\dfrac{-23}{p}\right)= -1$, then for all $\alpha\geq0$,
		\begin{align}\label{cong-A23-pnr}A_{23}\left(2p^{2\alpha+1}(pn+r)+ 2 p^{2(\alpha+1)}+1\right) \equiv 0 \pmod{2},		
	\end{align}
	where $r \in \{1,2,\ldots,p-1\}$.
	\end{theorem}

We organize the paper as follows. In the next section, we state a few well-known dissection formulas. In Section \ref{3}, we prove Theorems \ref{thm1.2}--\ref{thm1.4}. In Sections \ref{4}--\ref{7}, we prove Theorems \ref{thm1.6}--\ref{thm1.9}, respectively.
\section{Dissection formulas}\label{2}
Some known 2-, 3-, and 5-dissections formulas are stated in the following five lemmas.
 
\begin{lemma}\label{lem3} \cite[p. 40, Entries 25(i) and 25(ii)]{bcb3}
	If $\varphi$ is given by \eqref{varphi}, then
	\begin{align}\label{phi-2-dissect}
		\varphi(q) = \varphi(q^4) + 2q \psi (q^8).
	\end{align}
\end{lemma}
\begin{lemma}\cite[p. 315]{bcb3}\label{lem6}
	If $\psi$ is given by \eqref{psi}, then
	\begin{align}\label{psiqq3-2-dissect}
	\psi(q)\psi(q^3) &= \varphi(q^6)\psi(q^4) + q \varphi(q^2)\psi(q^{12}),\\
	\label{psiqq7-2-dissect}
		\psi(q)\psi(q^7) &= \varphi(q^{28})\psi(q^8) + q \psi(q^{2})\psi(q^{14}) + q^6 \varphi(q^{4})\psi(q^{56}).
	\end{align}
\end{lemma}

\begin{lemma}\label{lem8}\cite[p. 49 and p. 51]{bcb3}
	We have
	\begin{align}\label{psi3dissect}
		\psi(q) &= f(q^3,q^6) +q \psi(q^9),\\
		\label{varphi3dissect}\varphi(q)&=\varphi(q^9)+2q\psi(-q^9)\chi(q^3).
	\end{align}
\end{lemma}

\begin{lemma}\label{lem4}\cite[p. 49]{bcb3}
	We have
	\begin{align}\label{phi-5dissect}
		\varphi(q) = \varphi(q^{25}) + 2q f(q^{15}, q^{35}) + 2q^4 f(q^{5}, q^{45}).
	\end{align}
\end{lemma}

In the following lemma, we recall a $p$-dissection of $f_{1}$.
\begin{lemma}\label{lem1}\cite[Theorem 2.2]{CUI2013507}
	For a prime $p > 3$, we have 
	\begin{align*}
		f_{1} = (-1)^{\frac{\pm p-1}{6}} q^{\frac{p^2 - 1}{24}} f_{p^2} + \sum_{k \neq \frac{\pm p-1}{6}, k = - \frac{p-1}{2} }^{\frac{p-1}{2}} (-1)^k q^{\frac{3k^2 + k}{2}} f\Big(-q^{\frac{3p^2+(6k+1)p}{2}},-q^{\frac{3p^2-(6k+1)p}{2}}\Big),
	\end{align*} 
	where 
	\begin{center}
		$\dfrac{\pm p-1}{6}=$ 
		$\begin{cases}
			\dfrac{p-1}{6}, & if ~p\equiv~ 1~ (mod ~6),\\
			\dfrac{-p-1}{6}, & if ~p\equiv~ -1~ (mod ~6).
		\end{cases}$
	\end{center}
	Furthermore, for $\dfrac{-(p-1)}{k} \leq k \leq \dfrac{p-1}{k}$ and $k \neq \dfrac{\pm p-1}{6}$,
	\begin{align*}
		\dfrac{3k^2 + k}{2} \not\equiv \dfrac{p^2 - 1}{24} \pmod{p}.
	\end{align*}
\end{lemma}

\section{Proofs of Theorems \ref{thm1.2}--\ref{thm1.4}}\label{3}
\begin{proof}
With the aid of \eqref{chi}, we rewrite \eqref{genAn} in the form
	\begin{align}\label{equ1}
		\sum_{n=0}^{\infty} A(n) q^n = \chi(-q)\chi(-q^2).
	\end{align}
	
	From \cite[Lemma 3.5]{zbMATH02213594} (See also \cite[p. 194, Eq. (3.65)]{cooper}), we recall that
	\begin{align}\label{phiqq3}
		\varphi^2(q) - \varphi^2(q^3) = 4q \chi(q) \chi(-q^2)\psi(q^3)\psi(q^6).
	\end{align}
	Replacing $q$ by $-q$ in the above and then employing \eqref{varphi3dissect}, we have
	\begin{align}\label{equ2}
		4q \chi(-q) \chi(-q^2)\psi(-q^3)\psi(q^6)=\varphi^2(-q^3) -\left(\varphi(-q^9) -2q \psi(q^9) \chi(-q^3)\right)^2.
	\end{align}
	It follows from \eqref{equ1} and \eqref{equ2} that
	\begin{align}\label{gen-An}
		4\sum_{n=0}^{\infty} A(n) q^{n+1} = \frac{\varphi^2(-q^3)}{\psi(-q^3)\psi(q^6)} - \dfrac{1}{\psi(-q^3)\psi(q^6)}\left(\varphi(-q^9) -2q \psi(q^9) \chi(-q^3)\right)^2.
	\end{align}
Extracting the terms involving $q^{3n+1}$ from both sides of the above, we have		
\begin{align*}
		4\sum_{n=0}^{\infty} A(3n) q^{3n+1} = 4q\dfrac{\chi(-q^3)\psi(q^9)\varphi(-q^9)}{\psi(-q^3)\psi(q^6)}.
	\end{align*}
Dividing both sides of the above by $4q$ and then replacing $q^3$ by $q$, we arrive at
\begin{align*}
		\sum_{n=0}^{\infty} A(3n) q^n = \dfrac{\chi(-q)\psi(q^3)\varphi(-q^3)}{\psi(-q)\psi(q^2)},
	\end{align*}
	which is \eqref{gen-A3n}.
	
	Similarly, extracting the terms involving $q^{3n+2}$ from both sides of \eqref{gen-An}, dividing by $4q^2$ and then replacing $q^3$ by $q$, we arrive at \eqref{gen-A3n1}.
	
	Again, extracting the terms involving $q^{3n+3}$ from both sides of \eqref{gen-An}, we have
	\begin{align*}
		4\sum_{n=0}^{\infty} A(3n+2) q^{3n+3} = \dfrac{1}{\psi(-q^3)\psi(q^6)} \left(\varphi^2(-q^3)- \varphi^2(-q^9)\right).
	\end{align*}
	Replacing $q^3$ by $q$ in the above and then employing \eqref{phiqq3}, we find that
	\begin{align*}
		4\sum_{n=0}^{\infty} A(3n+2) q^{n+1} &=-4q\chi(-q)\chi(-q^2) \frac{\psi(-q^3)\psi(q^6)}{\psi(-q)\psi(q^2)},
	\end{align*}
	from which it follows that
	\begin{align*}
		\sum_{n=0}^{\infty} A(3n+2) q^n &=- \dfrac{\psi(-q^3)\psi(q^6)}{\psi^2(q^2)},
	\end{align*}
which is \eqref{gen-A3n2}.

Now, from \cite[Lemma 1.5]{da2022elementary}, we recall that
	\begin{align}\label{equ7}
		\dfrac{1}{\psi^2(q^2)} = \left(\frac{f_{6}^4f_{18}^6}{f_{12}^{12}} - 2q^6 \frac{f_{6}^7f_{36}^9}{f_{12}^{16}f_{18}^3}\right) + 3q^4\frac{f_{6}^6f_{36}^6}{f_{12}^{14}} - \left(2q^2\frac{f_{6}^5f_{18}^3f_{36}^3}{f_{12}^{13}} - q^8 \frac{f_{6}^8f_{36}^{12}}{f_{12}^{16}f_{18}^6} \right).
	\end{align}
	Employing \eqref{equ7} in \eqref{gen-A3n2}, extracting the terms involving $q^{3n+1}$, dividing both sides by $q$, and then replacing $q^3$ by $q$, we obtain
	\begin{align*}
		\sum_{n=1}^{\infty} A(9n+5)q^n &= -3q\psi(-q)\psi(q^2)\dfrac{f_{2}^6f_{12}^6}{f_{4}^{14}}= -3q\dfrac{f_1f_4}{f_2}\cdot\dfrac{f_4^2}{f_2}\cdot\dfrac{f_{2}^6f_{12}^6}{f_{4}^{14}}\\
		&=-3q\dfrac{f_1f_2^4f_{12}^6}{f_4^{11}},
	\end{align*}
	which is \eqref{gen-A9n5}.

Now, it is easy to see that
\begin{align*}
		\psi^3(q)\equiv \psi(q^3)\pmod3.
	\end{align*}
Employing the above in \eqref{gen-A3n2}, we find that
\begin{align}\label{eq-A3n2}
		\sum_{n=0}^{\infty} A(3n+2) q^n &=- \dfrac{\psi(-q^3)\psi(q^6)}{\psi^2(q^2)}=- \dfrac{\psi(-q^3)\psi(q^6)\psi(q^2)}{\psi^3(q^2)}\notag\\
		&\equiv 2\psi(-q^3)\psi(q^2) \pmod3,
	\end{align}
	which by \eqref{psi3dissect} can be written as
\begin{align*}
		\sum_{n=0}^{\infty} A(3n+2) q^n &\equiv 2\psi(-q^3)\left(f(q^6,q^{12}) +q^2 \psi(q^{18})\right) \pmod3.
	\end{align*}
	Extracting  the terms involving $q^{3n+2}$ from the above and then employing \eqref{psi3dissect} once again, we find that
\begin{align}\label{equ8}
		\sum_{n=0}^{\infty} A(9n+8) q^n &\equiv 2\psi(q^6)\psi(-q)\notag\\
		&\equiv 2\psi(q^6)\left(f(-q^3,q^{6}) -q \psi(-q^9)\right) \pmod3.
	\end{align}
Equating the coefficients of $q^{3n+2}$ from both sides of the above, we arrive at
	\begin{align*}
		A(27n+26) \equiv 0 \pmod3,
	\end{align*}
	which is \eqref{cong-A27n26}.
	
	Next, from \eqref{gen-A3n1}, \eqref{psi}, \eqref{chi} and the fact that $f_1^3\equiv f_3\pmod3$, we find that
	\begin{align}\label{A3n1-cong3}
		\sum_{n=0}^{\infty} A(3n+1) q^n&= -\dfrac{\psi^2(q^3)\chi^2(-q)}{\psi(-q)\psi(q^2)}=-\dfrac{f_1f_6^4}{f_3^2f_4^3}\notag\\		
		&\equiv 2 \frac{f_1f_6f_{18}}{f_3^2f_{12}} \pmod3.
	\end{align} 
	Again, extracting the terms involving $q^{3n}$ from both sides of \eqref{equ8} and then replacing $q^{3}$ by $q$, we arrive at
	\begin{align}\label{A27n8-cong3}
		\sum_{n=0}^{\infty} A(27n+8) q^n \equiv 2 \psi(q^2) f(-q,q^2) \pmod3.
	\end{align}

By \eqref{jtpi}, \eqref{varphi}, and \eqref{chi}, we find that
\begin{align*}
f(-q,q^2)=(q;-q^3)_\infty (-q^2;-q^3)_\infty(-q^3;-q^3)_\infty=\dfrac{(-q^3;q^6)_\infty^2(q^6;q^6)_\infty}{(-q;q^2)_\infty}.
\end{align*}
Employing the above identity and \eqref{psi} in \eqref{A27n8-cong3}, and then  simplifying by using the fact $f_1^3\equiv f_3\pmod3$ again, we find that
\begin{align}\label{A27n8-cong3a}
		\sum_{n=0}^{\infty} A(27n+8) q^n &\equiv 2 \dfrac{(q^4;q^4)_\infty(-q^3;q^6)_\infty^2(q^6;q^6)_\infty}{(-q;q^2)_\infty(q^2;q^4)_\infty}\notag\\
&\equiv 2 \dfrac{f_1f_4^3f_6^5}{f_2^3f_3^2f_{12}^2}\notag\\
&\equiv 2 \dfrac{f_1f_6f_{18}}{f_3^2f_{12}}		
		 \pmod3.
	\end{align}
From \eqref{A3n1-cong3} and \eqref{A27n8-cong3a}, we conculde that, for all $n\geq0$,
\begin{align*}A(3n+1)\equiv A(27n+8)		
		 \pmod3,
	\end{align*}	
which is \eqref{Silva-Sellers-Cong1}.

Next, extracting the terms involving $q^{3n+1}$ from both sides of \eqref{equ8}, dividing both sides by $q$ and then replacing $q^3$ by $q$, we find that
	\begin{align}\label{equ27n+17}
		\sum_{n=0}^{\infty} A(27n+17) q^n \equiv  \psi(q^2)\psi(-q^3) \pmod3.
	\end{align}	
	From the above congruence and \eqref{eq-A3n2}, we see that, for  all $n\geq0$, 
\begin{align*}
		A(3n+2)\equiv-A(27n+17)\pmod3,
	\end{align*}
	which is \eqref{cong-3n2-27n17}.

With the aid of \eqref{psi3dissect}, we can rewrite \eqref{equ27n+17} as 
\begin{align*}
		\sum_{n=0}^{\infty} A(27n+17) q^n &\equiv  \psi(-q^3)\left(f(q^6,q^{12}) +q^2 \psi(q^{18})\right) \pmod3.
	\end{align*}
	Equating the coefficients of $q^{3n+1}$ from both sides of the above, we find that
	\begin{align*}
		A(81n+44) \equiv 0 \pmod3,
	\end{align*}
	which is \eqref{cong-81n44}.
	
	Successive iterations of \eqref{cong-3n2-27n17}  give
\begin{align}\label{gencong-3n2}
		A(3n+2)&\equiv -A(3(9n+5)+2)\notag\\
		&\equiv A(27(9n+5)+17)\notag\\
		&\equiv A(3^5n + 3^3\cdot5 + 3\cdot5 +2)\notag\\
		&~\vdots\notag\\
		&\equiv (-1)^{j}A(3\cdot9^{j}n + 3\cdot9^{j-1}\cdot5 + 3.9^{j-2}\cdot5 +\cdots+ 3\cdot5 +2)\notag\\
		&\equiv (-1)^{j}A\left(3\cdot9^{j}n + \frac{15\cdot9^{j}+1}{8}\right)\pmod3.
	\end{align}
Replacing $n$ by $3n+1$ in the above, we find that
	\begin{align*}
		A(9n+5)\equiv (-1)^{j}A\left(9^{j+1}n + \dfrac{39\cdot9^j + 1}{8}\right) \pmod3.
	\end{align*}
Employing \eqref{cong-A9n5}, we see that, for all $j \geq 0$ and $n \geq 0$,
	\begin{align*}
		A\left(9^{j+1}n + \frac{39\cdot9^j + 1}{8}\right) \equiv 0 \pmod3,
	\end{align*}
	which is \eqref{Silva-Sellers-Cong2}.

Again, replacing $n$ by $9n+8$ in \eqref{gencong-3n2}, we have 
	\begin{align*}
		A(27n+26) \equiv (-1)^{j}A\left(3\cdot9^{j+1}n + \dfrac{23\cdot9^{j+1} + 1}{8}\right) \pmod3.
	\end{align*}
Employing \eqref{cong-A27n26} in the above, we readily arrive at \eqref{Silva-Sellers-Cong3}. 	
	
\end{proof}	
	
\section{Proof of Theorem \ref{thm1.6}}\label{4}
Setting $k=3$ in \eqref{gen-Akn} and then manipulating the $q$-products, we have
	\begin{align*}
		\sum_{n=0}^{\infty} A_{3}(n)q^n &= \chi(-q)\chi(-q^3) = \dfrac{f_{1}f_{3}}{f_{2}f_{6}}= \dfrac{\psi(-q)\psi(-q^3)}{f_4f_{12}}.
	\end{align*}
	Replacing $q$ by $-q$ in \eqref{psiqq3-2-dissect} and then using the resulting identity in the above, we have
	\begin{align*}
		\sum_{n=0}^{\infty} A_3(n)q^n & = \dfrac{\varphi(q^6)\psi(q^4)- q \varphi(q^2)\psi(q^{12})}{f_{4}f_{12}}. 
	\end{align*}
	Extracting, in turn, the even and odd terms from both sides of the above,  and then using \eqref{phi-2-dissect}, we find that
	\begin{align*}
		\sum_{n=0}^{\infty} A_{3}(2n)q^n &= \dfrac{\varphi(q^3)\psi(q^2)}{f_{2}f_{6}}= \dfrac{\psi(q^2)\left(\varphi(q^{12})+2q^3\psi(q^{24})\right)}{f_{2}f_{6}}\\\intertext{and}
		\sum_{n=0}^{\infty} A_{3}(2n+1)q^n &= - \frac{\varphi(q)\psi(q^6)}{f_{2}f_{6}}=- \dfrac{\psi(q^6)\left(\varphi(q^4)+2q\psi(q^8)\right)}{f_{2}f_{6}}.
\end{align*}
	Equating the coefficients of $q^{2n+1}$ from both sides of the above two identities, we arrive at
	\begin{align*}
	A_{3}(4n+r) \equiv 0 \pmod{2},\quad \textup{where}\quad r \in \{2, 3\},
	\end{align*}
	which is \eqref{cong-A3-4nr}.
	
\section{Proof of Theorem \ref{thm1.7}}\label{5}
Setting $k=5$ in \eqref{gen-Akn}, we have
\begin{align}\label{gen-A5n}
	\sum_{n=0}^{\infty}A_{5}(n)q^n &= \chi(-q)\chi(-q^5).
\end{align}

Now, recall from  \cite[p. 258, Entry 9(vii) and p. 262, Entry 10(iv)]{bcb3} that
\begin{align}\label{phiqq5}
	\varphi^2(q) - \varphi^2(q^5) &= 4qf(q,q^9)f(q^3,q^7)= 4q\chi(q) f_5f_{20}.
\end{align}
Multiplying by $\chi(q^5)$, and then replacing $q$ by $-q$, we find that
\begin{align*}
	4q\chi(-q)\chi(-q^5)&= \dfrac{\chi(-q^5)}{\chi(q^5)f_{10}f_{20}}\left(\varphi^2(-q^5)-\varphi^2(-q)\right)\\
	&=\dfrac{\varphi^2(-q^5)-\varphi^2(-q)}{\psi^2(q^5)}.
\end{align*}
With the aid of \eqref{gen-A5n} and \eqref{phi-5dissect}, the above may be rewritten as
\begin{align}
	&4\sum_{n=0}^{\infty}A_{5}(n)q^{n+1} \nonumber\\
	&=\frac{1}{\psi^{2}(q^5)}	\left(\varphi^2(-q^5) -  \left(\varphi(-q^{25}) - 2q f(-q^{15}, -q^{35}) + 2q^4 f(-q^{5}, -q^{45})\right)^2\right).\label{eq4}
\end{align}
Extracting, in turn, the terms of the form $q^{5n+2}$ and $q^{5n+3}$ from both sides of the above, we obtain
	\begin{align*}
		\sum_{n=0}^{\infty}A_{5}(5n+1)q^{n} &=  -\frac{f^2(-q^3,-q^7)}{\psi^{2}(q)}\equiv \frac{f(-q^{6},-q^{14})}{\psi(q^2)} \pmod{2},\\
	\sum_{n=0}^{\infty}A_{5}(5n+2)q^{n} &=  -q\frac{f^2(-q,-q^9)}{\psi^{2}(q)}\equiv q\frac{f(-q^{2},-q^{18})}{\psi(q^2)} \pmod{2},
	\end{align*}
	from which it readily follows that $A_{5}(10n+6) \equiv 0 \pmod2$ and $A_{5}(10n+2) \equiv 0 \pmod2$. This completes the proof of \eqref{cong-A5-10nr}.

Now, extracting the terms involving $q^{5n+5}$ from both sides of \eqref{eq4}, replacing $q^5$ by $q$, and then applying \eqref{phiqq5}, we find that
\begin{align*}
		4\sum_{n=0}^{\infty}A_{5}(5n+4)q^{n+1} &= \frac{1}{\psi^{2}(q)}	(\varphi^2(-q) -  \varphi^2(-q^{5})+ 8q f(-q^{3}, -q^{7}) f(-q, -q^{9}))\\
		&= \frac{4q f(-q^{3}, -q^{7}) f(-q, -q^{9})}{\psi^{2}(q)}\\
				&= \frac{4q\chi(-q)(-q^5;q^{10})_{\infty}f_{10}f_{20}}{\psi^{2}(q)}.\\
	\end{align*}
	Therefore,
	\begin{align}\label{gen-A5n4}
		\sum_{n=0}^{\infty}A_{5}(5n+4)q^{n} &= \dfrac{f_{1}^3f_{10}^3}{f_{2}^5f_{5}}\equiv \frac{f_{10}^2}{f_5}{f_{1}^3}\pmod5.
	\end{align}

But, well-known Jacobi's identity states that 	
	\begin{align*}
		f_{1}^3= \sum_{j=0}^{\infty} (-1)^j (2j+1) q^{{j(j+1)}/2}.
	\end{align*}
Employing this in \eqref{gen-A5n4}, we have 	
	\begin{align}\label{gen-A5n4mod5}
		\sum_{n=0}^{\infty}A_{5}(5n+4)q^{n} &\equiv \frac{f_{10}^2}{f_{5}} \sum_{j=0}^{\infty} (-1)^j (2j+1) q^{j(j+1)/2}\pmod{5}.
	\end{align}
	Now, $j(j+1)/2\equiv~0,1$ or $3\pmod 5$. Therefore, equating coefficients of $q^{5n+2}$ and $q^{5n+4}$, in turn, from both sides of \eqref{gen-A5n4mod5}, we find that	
\begin{align}\label{gen-A5n4mod5a}
		A_5(25n+14)\equiv A_{5}(25n+24)\equiv 0 \pmod5.
	\end{align}

Furthermore, if $j\equiv 2 \pmod5$, then $j(j+1)/2\equiv 3 \pmod5$ and $2j+1\equiv 0 \pmod5$. Therefore, equating the coefficients of $q^{5n+3}$ from both sides of \eqref{gen-A5n4mod5}, we find that
\begin{align}\label{gen-A5n4mod5b}
		A_5(25n+19)\equiv 0\pmod5.
	\end{align}
Clearly, \eqref{gen-A5n4mod5a} and \eqref{gen-A5n4mod5b} together give \eqref{cong-A5-25nr}. This completes the proof.

\section{Proof of Theorem \ref{thm1.8}}\label{6}

\noindent\emph{Proofs of \eqref{cong-A7-2n1-8n3}, \eqref{cong-A7-inf8n7}, \eqref{cong-A7-16nr}}. Setting $k=7$ in \eqref{gen-Akn}, manipulating the $q$-products, and then employing \eqref{psiqq7-2-dissect}, we have
\begin{align}\label{gen-A7n}
	\sum_{n=0}^{\infty}A_7(n)q^n &= \chi(-q)\chi(-q^7)=\dfrac{f_1f_7}{f_2f_{14}}=\dfrac{\psi(-q)\psi(-q^7)}{f_4f_{28}}\notag\\
	&=\dfrac{1}{f_4f_{28}}\left(\varphi(q^{28})\psi(q^8) - q \psi(q^2)\psi(q^{14}) + q^6 \varphi(q^4)\psi(q^{56})\right).
\end{align}
Extracting the odd terms from both sides of the above and then employing \eqref{psiqq7-2-dissect} once again, we have
\begin{align}\label{gen-A7-2n1}
	\sum_{n=0}^{\infty}A_7(2n+1)q^n &= -\dfrac{\psi(q)\psi(q^{7})}{f_2f_{14}}\notag\\
	&=-\dfrac{1}{f_2f_{14}}\left(\varphi(q^{28})\psi(q^8) + q \psi(q^2)\psi(q^{14}) + q^6 \varphi(q^4)\psi(q^{56})\right).
\end{align}
Extracting the odd terms, we find that
\begin{align}\label{gen-A7-4n3}
	\sum_{n=0}^{\infty}A_7(4n+3)q^n &= -\dfrac{\psi(q)\psi(q^{7})}{f_1f_{7}}=-\dfrac{f_2^2f_{14}^2}{f_1^2f_{7}^2}\notag\\
	&\equiv f_2f_{14}\pmod2.
\end{align}
It follows from \eqref{gen-A7-4n3} that
\begin{align}	\label{gen-A7-8n3}\sum_{n=0}^{\infty}A_7(8n+3)q^n &\equiv f_1f_{7}\pmod2\\\intertext{and}
\label{cong-A7-8n7}A_7(8n+7)	&\equiv 0\pmod2.
\end{align}

From \eqref{gen-A7-2n1}, we also have
\begin{align}\label{cong-A7-2n1}
	\sum_{n=0}^{\infty}A_7(2n+1)q^n &\equiv f_1f_{7}\pmod2.
\end{align}
From the above congruence and \eqref{gen-A7-8n3}, we readily arrive at \eqref{cong-A7-2n1-8n3}.

Now, iterating \eqref{cong-A7-2n1-8n3}, we find that
\begin{align*}
A_7(2n+1) &\equiv A_7(2(4n+1)+1)\\
&\equiv A_7(2(4^2n+4+1)+1)\\
&\equiv A_7(2(4^3n+4^2+4+1)+1)\\
&~\vdots\\
&\equiv A_7\left(2\left(4^jn+\dfrac{4^j-1}{3}\right)+1\right)\pmod2.
\end{align*}
Replacing $n$ by $4n+3$ in the above and then employing \eqref{cong-A7-8n7}, we obtain \eqref{cong-A7-inf8n7}.

Now, extracting the even terms on both sides of \eqref{gen-A7-2n1},  we have
\begin{align}\label{gen-A7-4n1}
	\sum_{n=0}^{\infty}A_7(4n+1)q^n &= -\dfrac{1}{f_1f_{7}}\left(\varphi(q^{14})\psi(q^4) + q^3 \varphi(q^2)\psi(q^{28})\right)\notag\\
	&=\dfrac{\psi(q)(\psi(q^7)\left(\varphi(q^{14})\psi(q^4) + q^3 \varphi(q^2)\psi(q^{28})\right)}{f_2^2f_{14}^2}.
	\end{align}

Now, as $\varphi(q)\equiv 1\pmod2$, from \eqref{psiqq7-2-dissect}, we have
\begin{align*}
\psi(q)\psi(q^7)
\equiv\psi(q^{8})+q\psi(q^2)\psi(q^{14})+q^6\psi(q^{56})\pmod2.
\end{align*}
Therefore, from\eqref{gen-A7-4n1}, we find that
\begin{align}\label{cong-A7-4n1}
	\sum_{n=0}^{\infty}A_7(4n+1)q^n &\equiv\dfrac{\left(\psi(q^{8})+q\psi(q^2)\psi(q^{14})+q^6\psi(q^{56})\right)\left(\psi(q^4) + q^3 \psi(q^{28})\right)}{f_4f_{28}}\pmod2.
	\end{align}
Extracting the even terms, we have
\begin{align*}
	&\sum_{n=0}^{\infty}A_7(8n+1)q^n \notag\\
	&\equiv\dfrac{\psi(q^2)\psi(q^4)+q^3\psi(q^2)\psi(q^{28})+q^2\psi(q)\psi(q^7)\psi(q^{14})}{f_2f_{14}}\notag\\
	&\equiv\dfrac{\psi(q^2)\psi(q^4)+q^3\psi(q^2)\psi(q^{28})+q^2\psi(q^{14})\left(\psi(q^{8})+q\psi(q^2)\psi(q^{14})+q^6\psi(q^{56})\right)}{f_2f_{14}}
	\pmod2.
	\end{align*}
	Extracting the odd terms from both sides of the above, we obtain
\begin{align*}
	\sum_{n=0}^{\infty}A_7(16n+9)q^n&\equiv\dfrac{q\psi(q)\psi(q^{14})+q\psi(q)\psi^2(q^7)}{f_1f_{7}}\notag\\
	&\equiv2\dfrac{q\psi(q)\psi(q^{14})}{f_1f_{7}}\equiv0
	\pmod2,
	\end{align*}
	from which \eqref{cong-A7-16nr} for $r=9$ follows readily.
	
	Next we prove \eqref{cong-A7-16nr} for $r=13$. Extracting the odd terms from both sides of \eqref{cong-A7-4n1}, we find that 
	\begin{align*}
	&\sum_{n=0}^{\infty}A_7(8n+5)q^n \notag\\
	&\equiv\dfrac{q\psi(q^4)\psi(q^{14})+q^4\psi(q^{14})\psi(q^{28})+\psi(q)\psi(q^2)\psi(q^7)}{f_2f_{14}}\notag\\
	&\equiv\dfrac{q\psi(q^4)\psi(q^{14})+q^4\psi(q^{14})\psi(q^{28})+
	\psi(q^2)\left(\psi(q^{8})+q\psi(q^2)\psi(q^{14})+q^6\psi(q^{56})\right)}{f_2f_{14}}\pmod2.
	\end{align*}
Extracting the odd terms, we obtain
\begin{align*}
	\sum_{n=0}^{\infty}A_7(16n+13)q^n &\equiv\dfrac{\psi(q^2)\psi(q^{7})+
	\psi^2(q)\psi(q^7)}{f_1f_7}\\
	&\equiv 2\dfrac{\psi(q^2)\psi(q^{7})}{f_1f_7}\equiv0\pmod2,
	\end{align*}
	from which \eqref{cong-A7-16nr} for $r=13$ is apparent. With this,  we complete the proof of \eqref{cong-A7-16nr}.

\noindent\emph{Proof of \eqref{cong-A7-14nr-inf}}. 	
At first, we show by the mathematical induction that for all $\alpha \geq 0$,
	\begin{align}\label{induct-state1}
		\sum_{n=0}^{\infty}	A_{7}\left(2\cdot7^{\alpha}n+\frac{2\cdot7^{\alpha}+1}{3}\right)q^n \equiv f_{1}f_{7} \pmod2.
	\end{align}
Clearly, by \eqref{cong-A7-2n1}, the result is true for $\alpha = 0$. Now, suppose that \eqref{induct-state1} holds good for some $\alpha>0$. Setting $p=7$ in Lemma \ref{lem1}, we have
	\begin{align*}
		f_{1} = q^{2} f_{49} + \sum_{k \neq 1, k = - 3 }^3 (-1)^k q^{\frac{3k^2 + k}{2}} f\left(-q^{\frac{3\cdot7^2+(6k+1)7}{2}},-q^{\frac{3\cdot7^2-(6k+1)7}{2}}\right).
	\end{align*}
	Employing the above in \eqref{induct-state1}, we have
\begin{align}\label{cong-A7-2n1-inf-7dissect}
	&\sum_{n=0}^{\infty}	A_{7}\left(2\cdot7^{\alpha}n+\frac{2\cdot7^{\alpha}+1}{3}\right)q^n \notag\\
	&\equiv q^2 f_7f_{49} + f_7\sum_{k \neq 1, k = - 3 }^3 (-1)^k q^{\frac{3k^2 + k}{2}} f\left(-q^{\frac{3\cdot7^2+(6k+1)7}{2}},-q^{\frac{3\cdot7^2-(6k+1)7}{2}}\right)\pmod2.
\end{align}
It can be easily verified that $\dfrac{3k^2+k}{2}\not\equiv2\pmod7$ for $k\neq1$. Therefore, extracting the terms involving $q^{7n+2}$ from both sides of the above, dividing both sides by $q^2$, and then replacing $q^7$ by $q$, we arrive at
	\begin{align*}
		\sum_{n=0}^{\infty} A_{7}\left(2\cdot7^{\alpha+1}n+\frac{2\cdot7^{\alpha+1}+1}{3}\right)q^n \equiv f_{1}f_{7} \pmod{2}.
	\end{align*}
Thus, \eqref{induct-state1} holds good for $\alpha+1$ whenever it holds good for some $\alpha>0$. Hence, by mathematical induction, \eqref{induct-state1} is true for all $\alpha\geq0$.

Now, it can also be seen that $(3k^2 + k)/2\equiv\ 0,1,2, \textup{or}~ 5 \pmod7$. Therefore, equating the coefficients of $q^{7n+r}$, where $
r=3,4,6$, from both sides of \eqref{cong-A7-2n1-inf-7dissect}, we arrive at
\begin{align*}A_7\left(2\cdot7^{\alpha}(7n+r)+\dfrac{2\cdot7^{\alpha}+1}{3}\right) &\equiv 0 \pmod2,\end{align*}
which is \eqref{cong-A7-14nr-inf}.

\noindent\emph{Proof of \eqref{cong-A7-pnr}}. First we prove by mathematical induction that 
if  $p$ is a prime such that $\left(\dfrac{-7}{p}\right) = -1$, then for all $\alpha \geq 0$ and $n \geq 0$,
	\begin{align}\label{cong-A7-pf1f7}
		\sum_{n=0}^{\infty}A_{7}\left(2\cdot p^{2\alpha}n+\dfrac{2\cdot p^{2\alpha}+1}{3}\right)q^n \equiv f_1f_7 \pmod2.
	\end{align}
	
	The case $\alpha=0$ of \eqref{cong-A7-pf1f7} is clearly true by \eqref{cong-A7-2n1}. 
	
	Suppose that \eqref{cong-A7-pf1f7} is true for some $\alpha>0$. Then, by Lemma \ref{lem1}, we have 
\begin{align}\label{cong-A7-pf1f7-a}
		&\sum_{n=0}^{\infty}A_{7}\left(2\cdot p^{2\alpha}n+\dfrac{2\cdot p^{2\alpha}+1}{3}\right)q^n \notag\\
		&\equiv \Big[\sum_{k \neq \frac{\pm p-1}{6}, k = - \frac{p-1}{2} }^{\frac{p-1}{2}} (-1)^k q^{\frac{3k^2 + k}{2}} f\Big(-q^{\frac{3p^2+(6k+1)p}{2}},-q^{\frac{3p^2-(6k+1)p}{2}}\Big)+(-1)^{\frac{\pm p-1}{6}} q^{\frac{p^2 - 1}{24}} f_{p^2}\Big] \notag\\
		&\quad\times\Big[\sum_{k \neq \frac{\pm p-1}{6}, k = - \frac{p-1}{2} }^{\frac{p-1}{2}} (-1)^k q^{7\cdot\frac{3k^2 + k}{2}} f\Big(-q^{7\cdot\frac{3p^2+(6k+1)p}{2}},-q^{7\cdot\frac{3p^2-(6k+1)p}{2}}\Big)\notag\\
		&\quad+(-1)^{\frac{\pm p-1}{6}} q^{7\cdot\frac{p^2 - 1}{24}} f_{7\cdot p^2} \Big]\pmod2.
	\end{align}
Now we consider the congruence 
	\begin{align}\label{cong-eqn}
		\frac{3k^2+k}{2}+7\cdot\frac{3m^2+m}{2} \equiv \frac{p^2-1}{3} \pmod{p},
	\end{align}
	where $-\frac{p-1}{2} \leq k , m \leq \frac{p-1}{2}$. Since the above congruence is equivalent to solving the congruence 
	\begin{align*}
		(6k+1)^2+7(6m+1)^2 \equiv 0 \pmod{p},
	\end{align*}
and $\left(\frac{-7}{p}\right) = -1$, it follows that \eqref{cong-eqn} has the unique solution $k = m = \frac{\pm p-1}{6}$. Therefore, extracting the terms involving $q^{pn+\frac{p^2-1}{3}}$ from both sides of \eqref{cong-A7-pf1f7-a}, we find that
	\begin{align}\label{cong-fp}
		\sum_{n=0}^{\infty} A_7\left(2\cdot p^{2\alpha+1}n+\frac{2 \cdot p^{2(\alpha+1)}+1}{3}\right)q^n \equiv f_pf_{7p} \pmod{2}.
	\end{align}
	Extracting the terms involving $q^{pn}$ from the above, we arrive at
	\begin{align*}
		\sum_{n=0}^{\infty} A_{7}\left(2\cdot p^{2(\alpha+1)}n+\frac{2\cdot p^{2(\alpha+1)}+1}{3}\right)q^n \equiv f_{1}f_{7} \pmod{2},
	\end{align*} 
	which clearly is the $\alpha+1$ case of \eqref{cong-A7-pf1f7}. This completes the proof of \eqref{cong-A7-pf1f7}.
	
	Now, equating the coefficients of $q^{pn+r}$ for $r\in\{1,2,\ldots,p-1\}$ on both sides of \eqref{cong-fp}, we readily arrive at \eqref{cong-A7-pnr}.

\begin{remark} It follows from \eqref{gen-A7n} and \eqref{gen-A7-2n1} that
	$$\left(\sum_{n=0}^{\infty} A_{7}(n)q^n	\right)\left(\sum_{n=0}^{\infty} A_{7}(2n+1)q^n\right) = -1.$$
\end{remark}

\section{Proof of Theorem \ref{thm1.9}}\label{7}
\noindent\emph{Proof of \eqref{cong-A23-23nr}}. At first, we prove by mathematical induction that for all $\alpha\geq0$, 
\begin{align}\label{cong-A23-gen-alpha}
		\sum_{n=0}^\infty A_{23}\left(2\cdot23^\alpha n+2\cdot23^\alpha+1\right)q^n \equiv f_1f_{23} \pmod2.
\end{align}

Setting $k=23$ in \eqref{gen-Akn}, we have
\begin{align}\label{gen-A23n}
	\sum_{n=0}^{\infty}A_{23}(n)q^n &= \chi(-q)\chi(-q^{23}).
\end{align}

From \cite[Eq. (7.4)]{zbMATH05199816}, we recall that
\begin{align*}
	\chi(-q)\chi(-q^{23})-\chi(q)\chi(q^{23})=-2q-2q^3(-q^2;q^2)_\infty(-q^{46};q^{46})_\infty,
\end{align*}
which, by \eqref{gen-A23n}, may be rewritten as
\begin{align*}
	\sum_{n=0}^{\infty}A_{23}(n)q^n -\sum_{n=0}^{\infty}A_{23}(n)(-q)^n&= -2q-2q^3(-q^2;q^2)_\infty(-q^{46};q^{46})_\infty.
\end{align*}
It follows from the above that 
\begin{align*}
	\sum_{n=0}^\infty A_{23}(2n+1)q^n &= -1-q(-q;q)_\infty(-q^{23};q^{23})_\infty\notag\\
	&\equiv 1+qf_1f_{23}\pmod2,
\end{align*}
and hence, 
\begin{align}\label{cong-A23-gen-alpha0}
	\sum_{n=0}^\infty A_{23}(2n+3)q^n &\equiv f_1f_{23}\pmod2,
\end{align}
which is the case $\alpha=0$ of \eqref{cong-A23-gen-alpha}.

Now, suppose that \eqref{cong-A23-gen-alpha} is true for some $\alpha>0$. We claim that it is then true for $\alpha+1$ as well. 

Setting $p=23$ in the $p$-dissection of $f_{1}$ stated in Lemma \ref{lem1}, we see that
	\begin{align*}
		f_{1} = q^{22} f_{23^2} + \sum_{k \neq -4, k = - 11 }^{11} (-1)^k q^{\frac{3k^2 + k}{2}} f\left(-q^{\frac{3\cdot 23^2+23(6k+1)}{2}},-q^{\frac{3\cdot23^2-23(6k+1)}{2}}\right).
	\end{align*}
	Employing the above in \eqref{cong-A23-gen-alpha}, we have
\begin{align}\label{cong-A23-gen-alpha-23}
		&\sum_{n=0}^\infty A_{23}\left(2\cdot23^\alpha n+2\cdot23^\alpha+1\right)q^n \notag\\
		&\equiv q^{22} f_{23}f_{23^2} + \sum_{k \neq -4, k = - 11 }^{11} (-1)^k q^{\frac{3k^2 + k}{2}} f_{23}f\left(-q^{\frac{3\cdot 23^2+23(6k+1)}{2}},-q^{\frac{3\cdot23^2-23(6k+1)}{2}}\right) \pmod2.
\end{align}
It is easy to verify that $\dfrac{3k^2 + k}{2}\not\equiv  22 \pmod{23}$ for $k\neq-4$. Therefore, extracting the terms involving $q^{23n+22}$ on both sides of the above, dividing by $q^{22}$, and then replacing $q^{23}$ by $q$, we find that
\begin{align*}
		\sum_{n=0}^\infty A_{23}\left(2\cdot23^\alpha (23n+22)+2\cdot23^\alpha+1\right)q^n &\equiv f_1f_{23}\pmod2,
\end{align*}
which is the $\alpha+1$ case of \eqref{cong-A23-gen-alpha}. Thus, \eqref{cong-A23-gen-alpha} holds good for all $\alpha\geq0$.

Now, we prove \eqref{cong-A23-23nr}. It can be easily verified that $\dfrac{3k^2 + k}{2}\not\equiv  4,6,9,10,13,14,16,18,19,20,21\pmod{23}$. So, equating the coefficients of $q^{23n+r}$ for $r\in\{4,6,9,10,13,14,16,18,19,20,21\}$ on both sides of \eqref{cong-A23-gen-alpha-23},  we find that, for all $\alpha\geq0$,
\begin{align*}
		A_{23}\left(2\cdot23^\alpha (23n+r)+2\cdot23^\alpha+1\right)\equiv 0\pmod2,
\end{align*}
which is \eqref{cong-A23-23nr}.

\noindent\emph{Proof of \eqref{cong-A23-pnr}}. We first prove by mathematical induction that if $p>3$ is a prime such that $\left(\dfrac{-23}{p}\right) = -1$, then for all $\alpha\geq0$
	\begin{align}\label{cong-f1f23-alpha}
		\sum_{n=0}^\infty A_{23}\left(2\cdot p^{2\alpha}n+2\cdot p^{2\alpha}+1\right)q^n \equiv f_1f_{23} \pmod2.
	\end{align}

Clearly, \eqref{cong-A23-gen-alpha0} is the $\alpha =0$ case of \eqref{cong-f1f23-alpha}.

Now, suppose that \eqref{cong-f1f23-alpha} is true for some $\alpha>0$. 
Then, by Lemma \ref{lem1}, we have 
\begin{align}\label{cong-A7-pf1f23-a}
		&\sum_{n=0}^{\infty}A_{23}\left(2\cdot p^{2\alpha}n+2\cdot p^{2\alpha}+1\right)q^n \notag\\
		&\equiv \Big[\sum_{k \neq \frac{\pm p-1}{6}, k = - \frac{p-1}{2} }^{\frac{p-1}{2}} (-1)^k q^{\frac{3k^2 + k}{2}} f\Big(-q^{\frac{3p^2+(6k+1)p}{2}},-q^{\frac{3p^2-(6k+1)p}{2}}\Big)+(-1)^{\frac{\pm p-1}{6}} q^{\frac{p^2 - 1}{24}} f_{p^2}\Big] \notag\\
		&\quad\times\Big[\sum_{k \neq \frac{\pm p-1}{6}, k = - \frac{p-1}{2} }^{\frac{p-1}{2}} (-1)^k q^{23\cdot\frac{3k^2 + k}{2}} f\Big(-q^{7\cdot\frac{3p^2+(6k+1)p}{2}},-q^{23\cdot\frac{3p^2-(6k+1)p}{2}}\Big)\notag\\
		&\quad+(-1)^{\frac{\pm p-1}{6}} q^{23\cdot\frac{p^2 - 1}{24}} f_{23\cdot p^2} \Big]\pmod2.
	\end{align}
	
	Now, consider the congruence
	\begin{align*}
		\dfrac{3k^2+k}{2} + 23\cdot\dfrac{3m^2+m}{2} \equiv p^2 - 1 \pmod{p},
	\end{align*}
	where $-\frac{p-1}{2} \leq k , m \leq \frac{p-1}{2}$.
	As the above congruence is equivalent to solving the congruence 
	\begin{align*}
		(6k+1)^2 + 23 (6m+1)^2 \equiv 0 \pmod{p},
	\end{align*}
and $\left(\dfrac{-23}{p}\right) = -1$, it has a unique solution, namely, $k = m = \frac{\pm p-1}{6}$. Therefore, extracting the terms involving $q^{pn+p^2-1}$ on both sides of the congruence \eqref{cong-A7-pf1f23-a}, dividing by $q^{p^2-1}$, and then replacing $q^p$ by $q$, we arrive at
	\begin{align}\label{cong-A23-fp}
		\sum_{n=0}^{\infty}A_{23}\left(2\cdot p^{2\alpha+1} n+ 2 p^{2(\alpha+1)}+1\right)q^n \equiv f_{p}f_{23p} \pmod{2}.
	\end{align}
	Extracting the terms involving $q^{pn}$ from both sides of the above and then replacing $q^p$ by $q$, we find that
	\begin{align*}
		\sum_{n=0}^{\infty}A_{23}\left(2\cdot p^{2(\alpha+1)}n + 2 p^{(2\alpha+1)}+1\right)q^n \equiv f_{1}f_{23} \pmod{2},
	\end{align*}
	which is clearly the $\alpha+1$ case of \eqref{cong-f1f23-alpha}. Hence, \eqref{cong-f1f23-alpha} is true for all $\alpha\geq0$.
	
Equating the coefficients of $q^{pn+r}$ for $r\in\{1,2,\ldots,p-1\}$ on both sides of \eqref{cong-A23-fp}, we readily arrive at \eqref{cong-A23-pnr} to complete the proof.

\section{Acknowledgement}
The second author was partially supported by the institutional fellowship for doctoral research from Tezpur University, Napaam, India. The author thanks the funding institution.

\bibliographystyle{plain}

\end{document}